\documentclass[11pt]{article}
\usepackage[a4paper,top=24mm,bottom=24mm,left=30mm,right=30mm]{geometry}
\usepackage{graphicx}
\usepackage[T1]{fontenc}
\usepackage[utf8]{inputenc}
\usepackage{libertinus}
\usepackage{microtype}
\usepackage{amsmath,amssymb,amsthm}
\usepackage{tikz}
\usepackage{xcolor}
\usepackage{enumitem}
\usepackage{booktabs}
\usepackage[hidelinks]{hyperref}
\usetikzlibrary{calc,arrows.meta}

\definecolor{ink}{HTML}{23211F}
\definecolor{soft}{HTML}{6E675E}
\definecolor{rule}{HTML}{C9C2B6}
\definecolor{domfill}{HTML}{EDE8DF}
\definecolor{domline}{HTML}{9A9184}
\definecolor{w2fill}{HTML}{F0D9A8}
\definecolor{w2line}{HTML}{B07C1E}
\definecolor{wrapfill}{HTML}{CFE2E4}
\definecolor{wrapline}{HTML}{2A7480}
\definecolor{carfillA}{HTML}{F6DED9}
\definecolor{carlineA}{HTML}{A8382C}
\definecolor{carfillB}{HTML}{DCE3EF}
\definecolor{carlineB}{HTML}{35508F}
\definecolor{shadec}{HTML}{F1EEE7}

\color{ink}

\tikzset{
  grid/.style     ={draw=rule,line width=.25pt},
  border/.style   ={draw=soft,line width=.7pt},
  dom/.style      ={fill=domfill,draw=domline,line width=.7pt,rounded corners=1.1pt},
  w2dom/.style    ={fill=w2fill,draw=w2line,line width=1pt,rounded corners=1.1pt},
  wrapdom/.style  ={fill=wrapfill,draw=wrapline,line width=1pt,rounded corners=1.1pt},
  wraparc/.style  ={draw=wrapline,line width=.7pt,densely dashed},
  hole/.style     ={draw=carlineA,fill=carfillA,line width=.8pt,densely dotted,rounded corners=1.1pt},
  flipmark/.style ={draw=soft,line width=.9pt,densely dashed,rounded corners=1.6pt},
  carfillA/.style ={fill=carfillA},
  carfillB/.style ={fill=carfillB},
  carpathA/.style ={draw=carlineA,line width=1.3pt,line cap=round},
  carpathB/.style ={draw=carlineB,line width=1.3pt,line cap=round},
  carstubA/.style ={carpathA,-{Stealth[length=3pt,width=3pt]}},
  carstubB/.style ={carpathB,-{Stealth[length=3pt,width=3pt]}},
  corner/.style   ={fill=w2line},
  lab/.style      ={text=soft},
  cellnote/.style ={text=ink,font=\scriptsize,inner sep=1pt,fill=white,fill opacity=.8,text opacity=1,rounded corners=1pt},
}

\usepackage{titlesec}
\titleformat{\section}{\large\bfseries}{\thesection.}{.6em}{}
\titleformat{\subsection}{\bfseries}{\thesubsection}{.6em}{}
\titlespacing*{\section}{0pt}{1.1\baselineskip}{.45\baselineskip}

\newtheoremstyle{plainish}{\topsep}{\topsep}{\itshape}{0pt}
  {\bfseries}{.}{.5em}{}
\theoremstyle{plainish}
\newtheorem{theorem}{Theorem}
\newtheorem{lemma}[theorem]{Lemma}

\theoremstyle{definition}

\theoremstyle{plainish}

\newcounter{fignum}
\renewcommand{\thefignum}{\arabic{fignum}}
\newcommand{\capt}[1]{{\footnotesize\color{soft}\refstepcounter{fignum}\textbf{Figure \thefignum.} #1}}
\newenvironment{figrow}{\par\smallskip\noindent\begin{center}}{\end{center}\par\smallskip}

\definecolor{cutcol}{HTML}{F3EFE8}
\definecolor{loopline}{HTML}{9A9184}

\tikzset{
  cutcol/.style    ={fill=cutcol},
  cutrule/.style   ={draw=w2line,line width=.8pt,densely dashed},
  labhi/.style     ={text=w2line,font=\bfseries},
pathA/.style     ={draw=carlineA,line width=1.3pt,line cap=round},
  pathAdot/.style  ={fill=carlineA},
  pathB/.style     ={draw=carlineB,line width=1.3pt,line cap=round},
  pathBdot/.style  ={fill=carlineB},
  loop/.style      ={draw=loopline,line width=1.1pt,line cap=round,densely dashed},
  loopdot/.style   ={fill=loopline},
crossodd/.style  ={fill=carfillA!55},
holeZ/.style     ={draw=domline,fill=none,line width=.6pt,densely dashed,rounded corners=1.1pt},
portZ/.style     ={draw=ink,fill=white,line width=.8pt},
  portR/.style     ={draw=ink,fill=ink,line width=.8pt},
rhoarc/.style    ={draw=wrapline,line width=1.1pt,line cap=round},
  rhoarrow/.style  ={rhoarc,-{Stealth[length=4pt,width=4pt]}},
  rhostub/.style   ={rhoarc,-{Stealth[length=4pt,width=4pt]}},
rhodom/.style    ={draw=ink,fill=none,line width=.7pt,rounded corners=1.1pt},
  ghostdom/.style  ={draw=rule,fill=white,line width=.5pt,rounded corners=1.1pt},
diag/.style      ={draw=soft,line width=.5pt,densely dotted},
  pathAlab/.style  ={text=carlineA},
  rholab/.style    ={text=wrapline},
  xout/.style      ={draw=carlineA,line width=1.4pt,line cap=round},
}

 \hypersetup{pdftitle={A Bijection for Domino Tilings of the M\"obius Band},
  pdfauthor={David Perkinson}}
\newtheorem{proposition}[theorem]{Proposition}
\newtheorem{corollary}[theorem]{Corollary}

\newcommand{\boardfig}[2]{\begin{minipage}[t]{\linewidth}\centering #1\\[4pt]\begin{minipage}{.84\linewidth}\capt{#2}\end{minipage}\end{minipage}}

\begin{document}

\begin{center}
{\LARGE\bfseries A Bijection for Domino Tilings of the M\"obius Band}\\[10pt]
{\large David Perkinson}\\[2pt]
{\color{soft}Reed College}\\[2pt]
{\color{soft}November 2026}\\[10pt]
\end{center}

\begin{center}
\begin{minipage}{.86\linewidth}
\small
\textbf{Abstract.}
A domino tiling of a rectangular board with an even number of rows and an even number
of columns is also a tiling of the M\"obius band obtained by gluing two opposite
borders of the board with a half-twist. The band has more tilings, because a domino
may cross the seam. Florescu, Morar, Perkinson, Salter, and Xu showed that the number
of tilings of the band equals a weighted count of the tilings of the rectangle, in
which the weight of each tiling is a power of two, and they asked for a bijective
proof. We give one. The tilings of the band are partitioned into classes, each class
contains exactly one tiling of the rectangle, and the size of the class is the weight
of that tiling. Combined with their evaluation of the weighted count, the bijection
proves the Lu--Wu product formula for the number of tilings of the band without
Pfaffians.

\medskip
{\footnotesize\color{soft}
}
\end{minipage}
\end{center}

\begin{center}
\includegraphics[width=.74\linewidth]{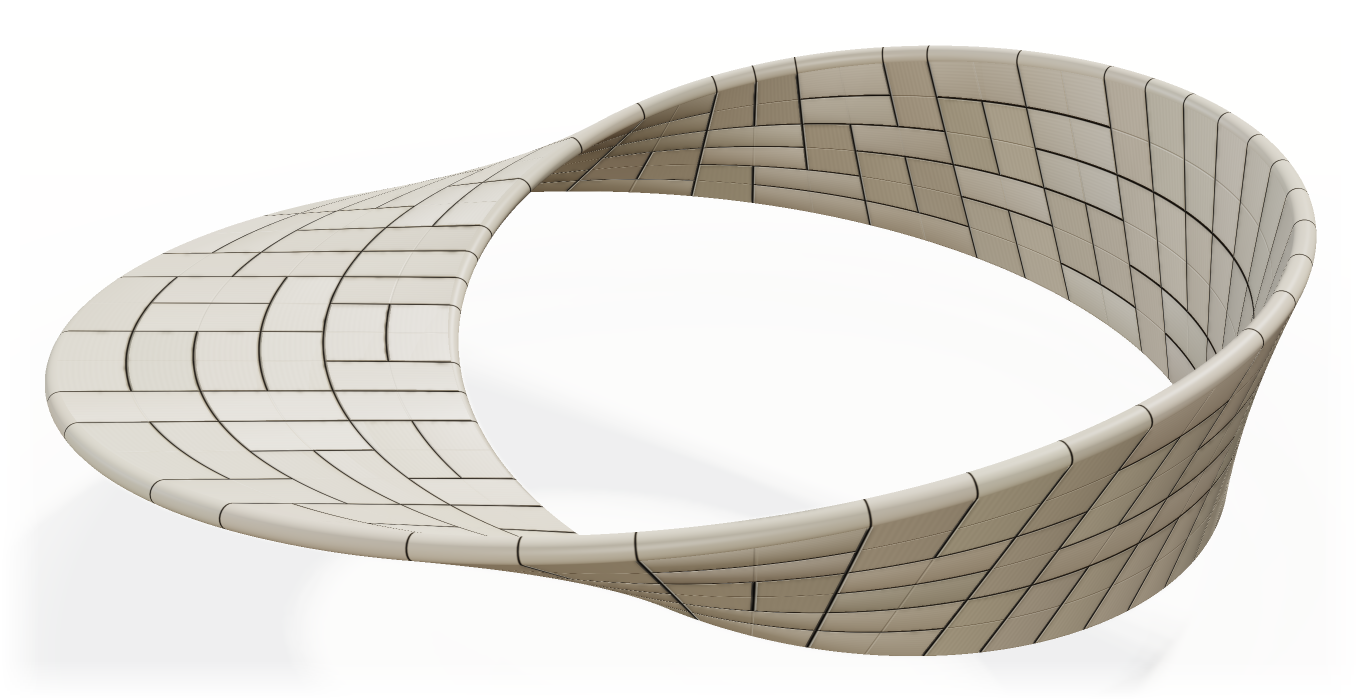}\\[4pt]
{\footnotesize\color{soft}A domino tiling of the M\"obius board with 6 rows and 48 columns.}
\end{center}

\vspace{2pt}

\section{Introduction}

Take the rectangular board~$R$ with $2m$ rows and $2n$ columns, and glue its left and
right borders with a half-twist: the right border of row~$r$ is glued to the left
border of row $\bar r := 2m+1-r$. The result is the \emph{M\"obius board}
$R^{\mathrm{mob}}$, a grid of $4mn$ square cells drawn on a M\"obius band. We call the identified
border the \emph{seam}.

\begin{figrow}
\boardfig{\begin{tikzpicture}[x=8mm,y=8mm,line join=round,line cap=round]

\draw[grid,step=1] (0,-6) grid (4,0);
\draw[border] (0,-6) rectangle (4,0);
\draw[border] (0,0) -- ++(-.14,0) (4,0) -- ++(.14,0);
\draw[border] (0,-1) -- ++(-.14,0) (4,-1) -- ++(.14,0);
\draw[border] (0,-2) -- ++(-.14,0) (4,-2) -- ++(.14,0);
\draw[border] (0,-3) -- ++(-.14,0) (4,-3) -- ++(.14,0);
\draw[border] (0,-4) -- ++(-.14,0) (4,-4) -- ++(.14,0);
\draw[border] (0,-5) -- ++(-.14,0) (4,-5) -- ++(.14,0);
\draw[border] (0,-6) -- ++(-.14,0) (4,-6) -- ++(.14,0);
\draw[border,-{Stealth[length=6.5pt,width=6pt]}] (0,-2.32) -- ++(0,-.36);
\draw[border,-{Stealth[length=6.5pt,width=6pt]}] (4,-3.68) -- ++(0,.36);
\draw[wrapdom] (3.13,-1.87) rectangle ++(0.74,0.74);
\draw[wrapdom] (0.13,-4.87) rectangle ++(0.74,0.74);
\node[text=ink,font=\footnotesize,text height=1.5ex,text depth=.25ex] at (-0.5,-0.5) {a};
\node[text=ink,font=\footnotesize,text height=1.5ex,text depth=.25ex] at (4.5,-5.5) {a};
\node[text=ink,font=\footnotesize,text height=1.5ex,text depth=.25ex] at (-0.5,-1.5) {b};
\node[text=ink,font=\footnotesize,text height=1.5ex,text depth=.25ex] at (4.5,-4.5) {b};
\node[text=ink,font=\footnotesize,text height=1.5ex,text depth=.25ex] at (-0.5,-2.5) {c};
\node[text=ink,font=\footnotesize,text height=1.5ex,text depth=.25ex] at (4.5,-3.5) {c};
\node[text=ink,font=\footnotesize,text height=1.5ex,text depth=.25ex] at (-0.5,-3.5) {d};
\node[text=ink,font=\footnotesize,text height=1.5ex,text depth=.25ex] at (4.5,-2.5) {d};
\node[text=ink,font=\footnotesize,text height=1.5ex,text depth=.25ex] at (-0.5,-4.5) {e};
\node[text=ink,font=\footnotesize,text height=1.5ex,text depth=.25ex] at (4.5,-1.5) {e};
\node[text=ink,font=\footnotesize,text height=1.5ex,text depth=.25ex] at (-0.5,-5.5) {f};
\node[text=ink,font=\footnotesize,text height=1.5ex,text depth=.25ex] at (4.5,-0.5) {f};
\node[lab,anchor=east] at (-1.05,-0.5) {\footnotesize 1};
\node[lab,anchor=east] at (-1.05,-1.5) {\footnotesize 2};
\node[lab,anchor=east] at (-1.05,-2.5) {\footnotesize 3};
\node[lab,anchor=east] at (-1.05,-3.5) {\footnotesize 4};
\node[lab,anchor=east] at (-1.05,-4.5) {\footnotesize 5};
\node[lab,anchor=east] at (-1.05,-5.5) {\footnotesize 6};
\node[lab,anchor=south] at (0.5,0.12) {\footnotesize 1};
\node[lab,anchor=south] at (1.5,0.12) {\footnotesize 2};
\node[lab,anchor=south] at (2.5,0.12) {\footnotesize 3};
\node[lab,anchor=south] at (3.5,0.12) {\footnotesize 4};
\end{tikzpicture}}{\label{fig:gluing}A $6\times4$ M\"obius board. Stretches of the border marked
with the same letter are glued. The letters run down the left border and up the right,
and this reversal is the half-twist. The two shaded cells are covered by a single
domino that crosses the seam.}
\end{figrow}

Every domino tiling of $R$ is also a tiling of $R^{\mathrm{mob}}$, but the M\"obius
board has more tilings, because a domino may now cross the seam. We regard the cells
as labeled, so two tilings that differ by a rotation of the band are counted
separately. How many tilings does the M\"obius board have?

One answer is a formula. A domino tiling is what a physicist calls a close-packed dimer
covering, and counting dimer coverings is one of the few exactly solvable problems of
two-dimensional statistical mechanics. Kasteleyn \cite{kasteleyn} and, independently,
Temperley and Fisher \cite{tf} solved the rectangle in 1961. Kasteleyn's method signs
the edges of the lattice so that the matchings are counted by a Pfaffian, the square
root of the determinant of an antisymmetric matrix. The signing depends on an
orientation of the surface, and a M\"obius band has none. What replaces the method on
nonorientable surfaces was settled in general by Tesler \cite{tesler}. The board
$R^{\mathrm{mob}}$ itself was first solved by Brankov and Priezzhev \cite{bp}. The
product formula for it is due to Lu and Wu \cite{luwu99,luwu}, and in the present
notation it reads
\begin{equation}\label{eq:luwu}
  \#\{\text{tilings of } R^{\mathrm{mob}}\} \;=\;
  \prod_{r=1}^{2m} \prod_{c=1}^{2n}
  \left[\, 4\cos^2\frac{r\pi}{2m+1} + 4\cos^2\frac{(2c-1)\pi}{4n} \,\right]^{1/4}.
\end{equation}
{\emergencystretch=1em
The 1961 formula for the rectangle~$R$ is the same product with $c\pi/(2n{+}1)$ in
place of $(2c{-}1)\pi/(4n)$. For the $4 \times 4$ board the products are $36$ for $R$
and $71$ for $R^{\mathrm{mob}}$.\par}

Another answer counts the tilings of $R^{\mathrm{mob}}$, the \emph{M\"obius tilings},
in terms of the tilings of $R$. We call a tiling of $R$ an \emph{ordinary} tiling. Let
$T$ be an ordinary tiling. Call a horizontal domino of $T$
occupying columns~$2n{-}1$ and $2n$ \emph{marked} if it lies in an even-numbered row,
and let $k(T)$ be the number of marked dominoes of $T$. Florescu, Morar, Perkinson,
Salter, and Xu \cite{sandpiles} proved that\footnote{For boards with two columns
($n = 1$), the count in \cite{sandpiles} is a weighted count of the tilings of the
$(2m{-}1) \times 2$ board. Identity~\eqref{eq:main} holds for these boards as written, and the proof below makes
no assumption on $n$.}
\begin{equation}\label{eq:main}
  \#\{\text{tilings of } R^{\mathrm{mob}}\} \;=\; \sum_{T} 2^{k(T)},
\end{equation}
the sum running over the tilings~$T$ of $R$. We call $2^{k(T)}$ the \emph{weight} of
$T$. For the $4 \times 4$ board the identity
reads
$11 \cdot 1 + 20 \cdot 2 + 5 \cdot 4 = 71$: eleven of the $36$ tilings of $R$ have
$k = 0$, twenty have $k = 1$, and five have $k = 2$.

The proof in \cite{sandpiles} evaluates the right side of \eqref{eq:main} through the
sandpile group of a grid graph, and observes that the result is the product
\eqref{eq:luwu}. So the identity appears there as a coincidence of two formulas, and
the authors ask for a bijection. They suggest a mechanism: pushing a marked domino
across the seam
``and then making room for this displacement by systematically shifting the other
dominos.''
The suggestion recalls the sliding-block puzzle Rush Hour, in which the cars in a
jammed lot are shuffled to let one car out. A tiling is a lot with no empty space at
all, so no domino can move alone. We will move the dominoes together, around a closed
loop.

This article gives the bijection. It matches the M\"obius tilings with the
pairs~$(T,M)$, where $T$ is an ordinary tiling and $M$ is a set of marked dominoes of
$T$. The only earlier bijective arguments for tilings of the band that we know of are
those of belcastro \cite{belcastro}, for boards two cells wide. These include the
M\"obius boards with $n = 1$.

The construction begins with a symmetry. The M\"obius board has a reflection about its
last column, defined by~\eqref{eq:tau} in Section~\ref{sec:companion}. The reflection
of a M\"obius tiling~$S$ is again a M\"obius tiling, the \emph{companion} $S^*$ of $S$.
Where $S$ has a horizontal domino occupying columns~$2n{-}1$ and $2n$, the companion instead
covers the same cell of the last column by a domino that crosses the seam, and
conversely. Figure~\ref{fig:running} shows a tiling and its companion.

The dominoes on which $S$ and $S^*$ differ lie along disjoint closed loops, and the
loops that cross the seam are the \emph{seam cycles} of $S$. Advancing the dominoes of
$S$ by one cell along a seam cycle gives another M\"obius tiling. We call this
operation a \emph{shift}, and we show that it does not change the seam cycles. So the
M\"obius tilings fall into \emph{shift classes}: the shift class of $S$ consists of the
tilings obtained from $S$ by shifting along some of its seam cycles, one tiling for
each set of seam cycles.

We prove that each shift class contains exactly one ordinary tiling~$T$, and that the
marked dominoes of $T$ lie on its seam cycles, exactly one on each. So the tilings in
the shift class of $T$ correspond to the sets of marked dominoes of $T$: shifting along
a seam cycle pushes its marked domino across the seam, and the other dominoes on the
cycle make room. The shift class of $T$ therefore has $2^{k(T)}$ members, and adding up
the sizes of the shift classes gives \eqref{eq:main}.

The bijection is Theorem~\ref{thm:phi}, and we reach it in three steps. After the
notation of Section~\ref{sec:prelim}, Section~\ref{sec:companion} constructs the reflection and the companion.
Section~\ref{sec:seam} proves that every seam cycle contains exactly two seam dominoes
(Proposition~\ref{thm:cars}), so that a shift changes the last column in exactly two rows.
Section~\ref{sec:bij} proves that every shift class contains exactly one ordinary
tiling (Proposition~\ref{thm:bij}), and deduces the bijection and identity~\eqref{eq:main}
(Corollary~\ref{cor:main}). The first two steps rest on
the reflection. The third adds a count of colors and the fact that the rectangle~$R$
is planar.
Section~\ref{sec:remarks} combines the
bijection with the evaluation of the right side of \eqref{eq:main} in \cite{sandpiles}
to prove the Lu--Wu formula \eqref{eq:luwu} without Pfaffians when the board has at
least four columns.

\section{Boards and matchings}\label{sec:prelim}

Rows are numbered from the top, columns from the left, and $(r,c)$ is the cell in row~$r$
and column~$c$. As
in the introduction, $\bar r = 2m+1-r$ is the row glued to row~$r$ at the seam. A
domino of $R$ covers two cells that share a side. On the M\"obius board a domino has $2m$
further positions, one for each row~$b$: the \emph{seam domino} $s_b$ covers $(b,2n)$ and
$(\bar b,1)$, the last cell of row~$b$ and the first cell of row~$\bar b$. The domino shaded in
Figure~\ref{fig:gluing} is $s_2$. A M\"obius tiling that uses no
seam domino is a tiling of $R$, that is, an ordinary tiling. A set of dominoes \emph{meets the seam} if it
contains a seam domino. We allow $n = 1$: nothing below assumes that the board has more than two
columns.

We record tilings as perfect matchings. Let $G^{\mathrm{mob}}$ be the graph whose
vertices are the cells and whose edges are the pairs of cells a single domino can
cover. A domino is an edge of $G^{\mathrm{mob}}$, and a M\"obius tiling is a perfect
matching. For tilings~$S$ and $S'$, the symmetric difference
$S \mathbin{\triangle} S'$ is the set of dominoes in exactly one of the two. The
dominoes on which two tilings differ lie along cycles: closed loops of cells on which
the dominoes of the two tilings alternate. On such a loop, advancing the dominoes of
one tiling by one cell turns them into the dominoes of the other. Doing this on any set
of the loops gives another tiling. These facts are standard. We state them as a lemma
and leave the straightforward proof to the reader.

\begin{lemma}\label{lem:super}
Let $S$ and $S'$ be perfect matchings of a graph. Then $S \mathbin{\triangle} S'$ is a
vertex-disjoint union of cycles, each alternating between $S$ and $S'$. For any set~$\mathcal{C}$ of these cycles, $S \mathbin{\triangle} \bigcup \mathcal{C}$ is again a
perfect matching.
\end{lemma}

\section{The reflection and the companion}\label{sec:companion}

The bijection moves dominoes around closed loops, and the loops come from comparing a
tiling with its mirror image. This section constructs the mirror image. We define a
reflection~$\tau$ of the M\"obius board, and with it the companion~$S^*$ of each
M\"obius tiling~$S$.

The last cell~$(b,2n)$ of row~$b$ can be covered in three ways: by a vertical domino
lying in the last column, by the seam domino~$s_b$, or by the horizontal domino that
also covers $(b,2n{-}1)$, which we call $h_b$. Fix a M\"obius tiling~$S$. Let
$B(S) \subseteq \{1, \dots, 2m\}$ be the set of rows~$b$ whose cell~$(b,2n)$ is
covered by $h_b$ or by $s_b$. We say that row~$b$ has \emph{type~$h$} or
\emph{type~$s$} accordingly, and we call a row of type~$s$ a \emph{seam row} of $S$.
The tiling~$S$ is ordinary exactly when it has no seam row. The number of rows in
$B(S)$ is even, because the other cells of the last column are covered in pairs by
vertical dominoes.

For an ordinary tiling~$T$, the marked dominoes of the introduction are the
dominoes~$h_b$ of $T$ with $b$ even, and $k(T)$ is their number. Our running example is
the ordinary tiling~$T$ of the $6 \times 4$ board on the left of
Figure~\ref{fig:running}. There the last column is covered by $h_1$ and $h_2$ and by
vertical dominoes in rows 3 to~6. So $B(T) = \{1,2\}$, and $k(T) = 1$ because $h_2$ is
marked.

Let $\Omega$ be the board with its last column removed. Define a map~$\tau$ on the cells
of the board by
\begin{equation}\label{eq:tau}
  \tau(r,c) \;=\; (\bar r,\,2n{-}c) \quad\text{for } c < 2n,
  \qquad\qquad \tau(r,\,2n) \;=\; (r,\,2n).
\end{equation}
The map~$\tau$ is an involution, and it fixes every cell of the last column. On $\Omega$
it is the half-turn of $\Omega$ about its center. It fixes no cell of $\Omega$, because
$\bar r \neq r$. So the cells fixed by $\tau$ are exactly the cells of the last column.

Seen on the band, $\tau$ is a reflection in the last column. For $1 \le j < 2n$, the
cell $j$ steps to the left of $(r,2n)$ along its row is $(r,\,2n{-}j)$, and the cell
$j$ steps to its right, across the seam, is $(\bar r,\,j)$. The map~$\tau$ exchanges
these two cells.

For a domino~$d$ we write $\tau(d)$ for the pair of images of its two cells, and we
apply $\tau$ to sets of dominoes.

\begin{lemma}\label{lem:tau}
The map~$\tau$ carries every domino to a domino. It carries each domino contained in
$\Omega$ to a domino contained in $\Omega$, it fixes each vertical domino of the last
column, and it exchanges $h_b$ with $s_b$ for every row~$b$.
\end{lemma}

\begin{proof}
A half-turn of $\Omega$ carries two cells of $\Omega$ that share a side to two cells of
$\Omega$ that share a side. A vertical domino of the last column is fixed because both
of its cells are. The domino~$h_b$ covers $(b,2n{-}1)$ and $(b,2n)$. Their images are
$(\bar b,1)$ and $(b,2n)$, which are the two cells of $s_b$. So $\tau(h_b) = s_b$, and
$\tau(s_b) = h_b$ because $\tau$ is an involution. Every domino is of one of these
kinds.
\end{proof}

The \emph{companion} of $S$ is its image under $\tau$:
\begin{equation}\label{eq:comp}
  S^* \;:=\; \tau(S) \;=\; \{\,\tau(d) : d \in S\,\}.
\end{equation}
By Lemma~\ref{lem:tau} the companion is a set of dominoes, and it covers every cell
exactly once because $\tau$ is a bijection on the cells. So the companion is a
M\"obius tiling. Since $\tau$ is an involution,
\[
  S^{**} \;=\; S.
\]
Lemma~\ref{lem:tau} also describes the companion in the last column. It has the same
vertical dominoes there as $S$. For each row~$b$ of $B(S)$, it contains $s_b$ if $S$
contains $h_b$, and it contains $h_b$ if $S$ contains $s_b$. So $B(S^*) = B(S)$, with
every type reversed. In particular, the companion of an ordinary tiling~$T$ contains
the seam domino~$s_b$ for every $b \in B(T)$.

In the running example, $\tau$ fixes column 4, reverses the rows of column 2, and
exchanges columns 1 and 3 with their rows reversed. The companion~$T^*$ is the middle
panel of Figure~\ref{fig:running}. It contains $s_1 = \tau(h_1)$, which covers $(1,4)$
and $(6,1)$, and $s_2 = \tau(h_2)$, which covers $(2,4)$ and $(5,1)$.

\begin{figrow}
\begin{minipage}[t]{\linewidth}\centering
\begin{tikzpicture}[x=8mm,y=8mm,line join=round,line cap=round]

\draw[grid,step=1] (0,-6) grid (4,0);
\draw[border] (0,-6) rectangle (4,0);
\draw[dom] (0.13,-0.87) rectangle ++(1.74,0.74);
\draw[dom] (2.13,-0.87) rectangle ++(1.74,0.74);
\draw[dom] (0.13,-1.87) rectangle ++(1.74,0.74);
\draw[w2dom] (2.13,-1.87) rectangle ++(1.74,0.74);
\draw[dom] (0.13,-3.87) rectangle ++(0.74,1.74);
\draw[dom] (1.13,-3.87) rectangle ++(0.74,1.74);
\draw[dom] (2.13,-3.87) rectangle ++(0.74,1.74);
\draw[dom] (3.13,-3.87) rectangle ++(0.74,1.74);
\draw[dom] (0.13,-4.87) rectangle ++(1.74,0.74);
\draw[dom] (2.13,-5.87) rectangle ++(0.74,1.74);
\draw[dom] (3.13,-5.87) rectangle ++(0.74,1.74);
\draw[dom] (0.13,-5.87) rectangle ++(1.74,0.74);
\node[lab,anchor=east] at (-0.35,-0.5) {\footnotesize 1};
\node[lab,anchor=east] at (-0.35,-1.5) {\footnotesize 2};
\node[lab,anchor=east] at (-0.35,-2.5) {\footnotesize 3};
\node[lab,anchor=east] at (-0.35,-3.5) {\footnotesize 4};
\node[lab,anchor=east] at (-0.35,-4.5) {\footnotesize 5};
\node[lab,anchor=east] at (-0.35,-5.5) {\footnotesize 6};
\node[lab,anchor=south] at (0.5,0.12) {\footnotesize 1};
\node[lab,anchor=south] at (1.5,0.12) {\footnotesize 2};
\node[lab,anchor=south] at (2.5,0.12) {\footnotesize 3};
\node[lab,anchor=south] at (3.5,0.12) {\footnotesize 4};
\end{tikzpicture}\hspace{1.6em}\begin{tikzpicture}[x=8mm,y=8mm,line join=round,line cap=round]

\draw[grid,step=1] (0,-6) grid (4,0);
\draw[border] (0,-6) rectangle (4,0);
\draw[dom] (0.13,-1.87) rectangle ++(0.74,1.74);
\draw[dom] (1.13,-0.87) rectangle ++(1.74,0.74);
\draw[wrapdom] (3.13,-0.87) rectangle ++(0.74,0.74);
\draw[wrapdom] (0.13,-5.87) rectangle ++(0.74,0.74);
\draw[dom] (1.13,-1.87) rectangle ++(1.74,0.74);
\draw[wrapdom] (3.13,-1.87) rectangle ++(0.74,0.74);
\draw[wrapdom] (0.13,-4.87) rectangle ++(0.74,0.74);
\draw[dom] (0.13,-3.87) rectangle ++(0.74,1.74);
\draw[dom] (1.13,-3.87) rectangle ++(0.74,1.74);
\draw[dom] (2.13,-3.87) rectangle ++(0.74,1.74);
\draw[dom] (3.13,-3.87) rectangle ++(0.74,1.74);
\draw[dom] (1.13,-4.87) rectangle ++(1.74,0.74);
\draw[dom] (3.13,-5.87) rectangle ++(0.74,1.74);
\draw[dom] (1.13,-5.87) rectangle ++(1.74,0.74);
\node[lab,anchor=east] at (-0.35,-0.5) {\footnotesize 1};
\node[lab,anchor=east] at (-0.35,-1.5) {\footnotesize 2};
\node[lab,anchor=east] at (-0.35,-2.5) {\footnotesize 3};
\node[lab,anchor=east] at (-0.35,-3.5) {\footnotesize 4};
\node[lab,anchor=east] at (-0.35,-4.5) {\footnotesize 5};
\node[lab,anchor=east] at (-0.35,-5.5) {\footnotesize 6};
\node[lab,anchor=south] at (0.5,0.12) {\footnotesize 1};
\node[lab,anchor=south] at (1.5,0.12) {\footnotesize 2};
\node[lab,anchor=south] at (2.5,0.12) {\footnotesize 3};
\node[lab,anchor=south] at (3.5,0.12) {\footnotesize 4};
\node[cellnote] at (3.5,-0.5) {$s_{1}$};
\node[cellnote] at (0.5,-5.5) {$s_{1}$};
\node[cellnote] at (3.5,-1.5) {$s_{2}$};
\node[cellnote] at (0.5,-4.5) {$s_{2}$};
\end{tikzpicture}\hspace{1.6em}\begin{tikzpicture}[x=8mm,y=8mm,line join=round,line cap=round]

\draw[grid,step=1] (0,-6) grid (4,0);
\draw[border] (0,-6) rectangle (4,0);
\fill[carfillA] (0,-1) rectangle ++(1,1);
\fill[carfillA] (1,-1) rectangle ++(1,1);
\fill[carfillA] (2,-1) rectangle ++(1,1);
\fill[carfillA] (3,-1) rectangle ++(1,1);
\fill[carfillA] (0,-2) rectangle ++(1,1);
\fill[carfillA] (1,-2) rectangle ++(1,1);
\fill[carfillA] (2,-2) rectangle ++(1,1);
\fill[carfillA] (3,-2) rectangle ++(1,1);
\fill[carfillA] (0,-5) rectangle ++(1,1);
\fill[carfillA] (1,-5) rectangle ++(1,1);
\fill[carfillA] (2,-5) rectangle ++(1,1);
\fill[carfillA] (0,-6) rectangle ++(1,1);
\fill[carfillA] (1,-6) rectangle ++(1,1);
\fill[carfillA] (2,-6) rectangle ++(1,1);
\draw[dom] (0.13,-0.87) rectangle ++(1.74,0.74);
\draw[dom] (2.13,-0.87) rectangle ++(1.74,0.74);
\draw[dom] (0.13,-1.87) rectangle ++(1.74,0.74);
\draw[w2dom] (2.13,-1.87) rectangle ++(1.74,0.74);
\draw[dom] (0.13,-3.87) rectangle ++(0.74,1.74);
\draw[dom] (1.13,-3.87) rectangle ++(0.74,1.74);
\draw[dom] (2.13,-3.87) rectangle ++(0.74,1.74);
\draw[dom] (3.13,-3.87) rectangle ++(0.74,1.74);
\draw[dom] (0.13,-4.87) rectangle ++(1.74,0.74);
\draw[dom] (2.13,-5.87) rectangle ++(0.74,1.74);
\draw[dom] (3.13,-5.87) rectangle ++(0.74,1.74);
\draw[dom] (0.13,-5.87) rectangle ++(1.74,0.74);
\draw[carpathA] (0.5,-0.5) -- (1.5,-0.5);
\draw[carpathA] (1.5,-0.5) -- (2.5,-0.5);
\draw[carpathA] (2.5,-0.5) -- (3.5,-0.5);
\draw[carstubA] (3.5,-0.5) -- (4.5,-0.5);
\draw[carstubA] (0.5,-5.5) -- (-0.5,-5.5);
\draw[carpathA] (0.5,-5.5) -- (1.5,-5.5);
\draw[carpathA] (1.5,-5.5) -- (2.5,-5.5);
\draw[carpathA] (2.5,-5.5) -- (2.5,-4.5);
\draw[carpathA] (2.5,-4.5) -- (1.5,-4.5);
\draw[carpathA] (1.5,-4.5) -- (0.5,-4.5);
\draw[carstubA] (0.5,-4.5) -- (-0.5,-4.5);
\draw[carstubA] (3.5,-1.5) -- (4.5,-1.5);
\draw[carpathA] (3.5,-1.5) -- (2.5,-1.5);
\draw[carpathA] (2.5,-1.5) -- (1.5,-1.5);
\draw[carpathA] (1.5,-1.5) -- (0.5,-1.5);
\draw[carpathA] (0.5,-1.5) -- (0.5,-0.5);
\fill[carlineA] (0.5,-0.5) circle (1.5pt);
\fill[carlineA] (1.5,-0.5) circle (1.5pt);
\fill[carlineA] (1.5,-0.5) circle (1.5pt);
\fill[carlineA] (2.5,-0.5) circle (1.5pt);
\fill[carlineA] (2.5,-0.5) circle (1.5pt);
\fill[carlineA] (3.5,-0.5) circle (1.5pt);
\fill[carlineA] (3.5,-0.5) circle (1.5pt);
\fill[carlineA] (0.5,-5.5) circle (1.5pt);
\fill[carlineA] (0.5,-5.5) circle (1.5pt);
\fill[carlineA] (1.5,-5.5) circle (1.5pt);
\fill[carlineA] (1.5,-5.5) circle (1.5pt);
\fill[carlineA] (2.5,-5.5) circle (1.5pt);
\fill[carlineA] (2.5,-5.5) circle (1.5pt);
\fill[carlineA] (2.5,-4.5) circle (1.5pt);
\fill[carlineA] (2.5,-4.5) circle (1.5pt);
\fill[carlineA] (1.5,-4.5) circle (1.5pt);
\fill[carlineA] (1.5,-4.5) circle (1.5pt);
\fill[carlineA] (0.5,-4.5) circle (1.5pt);
\fill[carlineA] (0.5,-4.5) circle (1.5pt);
\fill[carlineA] (3.5,-1.5) circle (1.5pt);
\fill[carlineA] (3.5,-1.5) circle (1.5pt);
\fill[carlineA] (2.5,-1.5) circle (1.5pt);
\fill[carlineA] (2.5,-1.5) circle (1.5pt);
\fill[carlineA] (1.5,-1.5) circle (1.5pt);
\fill[carlineA] (1.5,-1.5) circle (1.5pt);
\fill[carlineA] (0.5,-1.5) circle (1.5pt);
\fill[carlineA] (0.5,-1.5) circle (1.5pt);
\fill[carlineA] (0.5,-0.5) circle (1.5pt);
\node[lab,anchor=east] at (-0.95,-0.5) {\footnotesize 1};
\node[lab,anchor=east] at (-0.95,-1.5) {\footnotesize 2};
\node[lab,anchor=east] at (-0.95,-2.5) {\footnotesize 3};
\node[lab,anchor=east] at (-0.95,-3.5) {\footnotesize 4};
\node[lab,anchor=east] at (-0.95,-4.5) {\footnotesize 5};
\node[lab,anchor=east] at (-0.95,-5.5) {\footnotesize 6};
\end{tikzpicture}\\[6pt]
\begin{minipage}{.92\linewidth}\capt{The running example.\label{fig:running}
The left panel shows an ordinary tiling~$T$ of the $6\times4$ board, with its marked
domino~$h_2$ in gold. The middle panel shows the companion~$T^*$, in which each seam
domino appears as two labeled halves. The right panel shows in red the cycle along
which $T$ and $T^*$ differ, drawn on $T$.}
\end{minipage}\end{minipage}
\end{figrow}

\section{Seam cycles}\label{sec:seam}

Every M\"obius tiling~$S$ now has a companion~$S^*$. This section defines the seam
cycles of $S$, the loops around which the bijection moves dominoes. Its main result,
Proposition~\ref{thm:cars}, is that a seam cycle contains exactly two seam dominoes. So
moving the dominoes around a seam cycle changes the last column in exactly two rows.
Throughout, $S$ is a M\"obius tiling.

A \emph{seam cycle} of $S$ is a component of $S \mathbin{\triangle} S^*$ that contains
a seam domino. Since $S^{**} = S$, the tilings~$S$ and $S^*$ have the same seam
cycles. In the running example, $T$ and $T^*$ differ along a single cycle, and it is a
seam cycle (Figure~\ref{fig:running}, right).

Let $C$ be a seam cycle of $S$. By Lemma~\ref{lem:super}, replacing the dominoes of
$S$ on $C$ by the dominoes of $S^*$ on $C$ converts $S$ into another M\"obius tiling
$S \mathbin{\triangle} C$. We call this operation \emph{shifting along} $C$. Each
domino on the cycle advances by one cell. If $C$ contains the seam domino~$s_b$, then
it also contains $h_b$, because one of the two covers the cell~$(b,2n)$ in $S$ and the
other covers it in $S^*$. The shift exchanges them, and so it reverses the type of
row~$b$. It does not change $B(S)$ or the vertical dominoes of the last column. In
the running example, $C$ is the only component of
$T \mathbin{\triangle} T^*$, so $T \mathbin{\triangle} C = T^*$. Advancing the marked
domino~$h_2$ onto $s_2$ pushes the horizontal domino of row 6 onto $s_1$ and draws
$h_1$ back into $\Omega$. In general, $S \mathbin{\triangle} S^*$ has other components,
and shifting along $C$ does not give the companion. Figure~\ref{fig:D} shows a tiling
with two seam cycles and two components that do not meet the seam.

\begin{figrow}
\boardfig{\input{figs/Dcross.tex}}{The two seam cycles of an ordinary tiling~$T$ of the
$10 \times 20$ board with $B(T) = \{3,6,9,10\}$. The tiling itself is not drawn.
The rows of $B(T)$ are numbered in
gold. Each seam domino is drawn as two segments that leave the board, labeled with its
name. The dashed gray cycles are the components of
$T \mathbin{\triangle} T^*$ that do not meet the seam. The shading is used in
Section~\ref{sec:bij}.\label{fig:D}}
\end{figrow}

\begin{proposition}[Seam cycles]\label{thm:cars}
Let $S$ be a M\"obius tiling. Each seam cycle of $S$ contains exactly two seam
dominoes, and $S$ has exactly $|B(S)|/2$ seam cycles. If $C$ is a seam cycle of $S$,
then $S \mathbin{\triangle} C$ has the same seam cycles as $S$.
\end{proposition}

For the proof, we show that $\tau$ carries each seam cycle to itself and acts on it as
a reflection of a cycle graph of even length. Such a reflection fixes two cells or
none, and the proposition follows by counting the cells that $\tau$ fixes.

\begin{lemma}\label{lem:inv}
The map~$\tau$ carries $S \mathbin{\triangle} S^*$ to itself, and it carries each seam
cycle of $S$ to itself.
\end{lemma}

\begin{proof}
By Lemma~\ref{lem:tau}, $\tau$ is a bijection on the dominoes, and taking images under
a bijection commutes with symmetric difference. Hence
$\tau(S \mathbin{\triangle} S^*) = \tau(S) \mathbin{\triangle} \tau(S^*)
 = S^* \mathbin{\triangle} S$. Two dominoes share a cell exactly when their images do,
so $\tau$ carries each component of $S \mathbin{\triangle} S^*$ to a component. Let $C$
be a seam cycle. It contains a seam domino, so it passes through a cell~$x$ of the
last column. Since $\tau$ fixes $x$, the component~$\tau(C)$ also passes through $x$.
Distinct components are disjoint, so $\tau(C) = C$.
\end{proof}

\begin{lemma}\label{lem:two}
Let $C$ be a seam cycle of $S$. List the cells of $C$ in order round the cycle as
$x_0, x_1, \dots, x_{L-1}$, with $x_0$ in the last column and with subscripts read
modulo $L$. Then $L$ is even, and $\tau(x_i) = x_{-i}$ for every $i$. The cells of $C$
in the last column are $x_0$ and $x_{L/2}$. If these are the cells~$(b,2n)$ and
$(b',2n)$, then the dominoes of $C$ that cover them are $h_b$ and $s_b$ and $h_{b'}$ and
$s_{b'}$.
\end{lemma}

\begin{proof}
By Lemma~\ref{lem:super}, $C$ is a cycle whose dominoes alternate between $S$ and
$S^*$. So the cells $x_0, \dots, x_{L-1}$ are distinct, and each of them is covered by
exactly two dominoes of $C$, one of $S$ and one of $S^*$. The dominoes of $S$ on $C$
therefore cover each cell of $C$ exactly once. Since each domino covers two cells,
$L$ is even.

Let $x$ be a cell of $C$ in the last column, and let $a$ be its row. The two dominoes
of $C$ that cover $x$ are a domino~$d$ of $S$ and a domino of $S^*$. Since $\tau$ fixes
$x$, the domino of $S^*$ that covers $x$ is $\tau(d)$. These two dominoes are
different, because a domino of $C$ belongs to only one of $S$ and $S^*$. The domino~$d$ covers a cell of the last column in row~$a$, so it is a
vertical domino of the last column or $h_a$ or $s_a$. By Lemma~\ref{lem:tau} it is not
vertical, since $\tau(d) \neq d$. Hence $d$ and $\tau(d)$ are $h_a$ and $s_a$. This proves the last assertion of the lemma. It also
shows that $\tau$ exchanges the other cell of $h_a$ with the other cell of $s_a$. These
are the two cells next to $x$ on $C$.

In particular $\tau$ fixes $x_0$ and exchanges $x_1$ with $x_{-1}$. Call two cells of
$C$ \emph{adjacent in} $C$ if a domino of $C$ covers both, so that the cells adjacent
to $x_i$ in $C$ are $x_{i-1}$ and $x_{i+1}$. By Lemma~\ref{lem:inv}, $\tau$ carries
each domino of $C$ to a domino of $C$. So $\tau$ permutes the cells of $C$, and it
carries cells adjacent in $C$ to cells adjacent in $C$. Hence, for each $i$, the
cell~$\tau(x_{i+1})$ is adjacent to $\tau(x_i)$ in $C$ and is different from
$\tau(x_{i-1})$. Since $\tau(x_0) = x_0$ and $\tau(x_1) = x_{-1}$, it follows by
induction on $i$ that $\tau(x_i) = x_{-i}$ for every $i$.

Thus $\tau$ restricted to $C$ is a reflection of a cycle graph of even length. It
fixes exactly two cells of $C$, namely $x_0$ and the opposite cell~$x_{L/2}$, since
$i \equiv -i \pmod L$ only for $i = 0$ and $i = L/2$. The cells fixed by $\tau$ are
the cells of the last column. Hence
the cells of $C$ in the last column are $x_0$ and $x_{L/2}$.
\end{proof}

In the running example the seam cycle has $14$ cells (Figure~\ref{fig:running}, right).
Starting at $x_0 = (1,4)$ and going left along row~1, it returns to the last column at
$x_7 = (2,4)$. The map~$\tau$ exchanges $x_1, \dots, x_6$, which lie in rows 1 and 2,
with $x_{13}, \dots, x_8$, which lie in rows 6 and 5.

\begin{proof}[Proof of Proposition~\ref{thm:cars}]
Let $C$ be a seam cycle of $S$. By Lemma~\ref{lem:two}, the cycle~$C$ passes through
exactly two cells of the last column, say $(b,2n)$ and $(b',2n)$, and the dominoes of
$C$ that cover them are $h_b$ and $s_b$ and $h_{b'}$ and $s_{b'}$. Every seam domino
covers a cell of the last column. So $C$ contains exactly two seam dominoes, $s_b$ and
$s_{b'}$.

Next we count the seam cycles. For $b \in B(S)$, the cell~$(b,2n)$ is covered by $h_b$
in one of $S$ and $S^*$ and by $s_b$ in the other, so it lies on a seam cycle. A cell
of the last column whose row is not in $B(S)$ is covered by the same vertical domino in
$S$ and in $S^*$, so it lies on no cycle. The seam cycles are disjoint, and each passes
through two of the $|B(S)|$ cells~$(b,2n)$ with $b \in B(S)$. So there are $|B(S)|/2$
seam cycles.

Finally, let $C$ be a seam cycle of $S$. Since $\tau(C) = C$ by Lemma~\ref{lem:inv},
\[
  (S \mathbin{\triangle} C) \mathbin{\triangle} (S \mathbin{\triangle} C)^*
  \;=\; (S \mathbin{\triangle} C) \mathbin{\triangle} (\tau(S) \mathbin{\triangle} \tau(C))
  \;=\; (S \mathbin{\triangle} C) \mathbin{\triangle} (S^* \mathbin{\triangle} C)
  \;=\; S \mathbin{\triangle} S^*.
\]
So the cycles along which $S \mathbin{\triangle} C$ differs from its companion are the
cycles along which $S$ differs from $S^*$. In particular, $S \mathbin{\triangle} C$ and
$S$ have the same seam cycles.
\end{proof}

By Proposition~\ref{thm:cars}, a seam cycle of $S$ contains exactly two seam dominoes, say
$s_b$ and $s_{b'}$. We say that the seam cycle \emph{pairs} rows~$b$ and $b'$. It may
run through many rows of the board, but shifting along it reverses the types of
rows~$b$ and $b'$ and of no other row. For each $b \in B(S)$, the seam domino~$s_b$
belongs to exactly one of $S$ and $S^*$, so it lies on a seam cycle. Since the seam
cycles are disjoint, they divide the rows of $B(S)$ into pairs.

\section{Shift classes and the bijection}\label{sec:bij}

Shifting along seam cycles leads from one M\"obius tiling to others. This section shows
that every M\"obius tiling is reached in this way from exactly one ordinary tiling
(Proposition~\ref{thm:bij}). The bijection follows (Theorem~\ref{thm:phi}). Throughout,
$S$ is a M\"obius tiling.

Let $\mathcal{C}$ be a set of seam cycles of $S$. By the last assertion of
Proposition~\ref{thm:cars}, after a shift along one cycle of $\mathcal{C}$ the others are
still seam cycles, so we may shift along the cycles of $\mathcal{C}$ one after another.
The cycles are disjoint, so the result is $S \mathbin{\triangle} \bigcup \mathcal{C}$
whatever the order. The tilings so obtained form
the \emph{shift class} of $S$. Distinct sets~$\mathcal{C}$ reverse the types of
different sets of rows, so the shift class has $2^k$ members, where $k$ is the
number of seam cycles of $S$. A tiling in the shift class of $S$ has the same seam
cycles as $S$, and hence the same shift class. So distinct shift classes are disjoint, and the shift
classes partition the set of M\"obius tilings.

The components of $S \mathbin{\triangle} S^*$ that do not meet the seam play no part.
Exchanging dominoes along one of them leaves the last column alone and produces a
tiling in a different shift class.

An ordinary tiling is a tiling with no seam row. To find one in the shift class of $S$,
we shift along seam cycles until no seam row is left. For this the two rows paired by a
seam cycle must have the same type. If a seam cycle paired a seam row with a row of
type~$h$, shifting along it would only exchange the two types, and every tiling in the shift
class would have a seam row. The proof that the two types agree
(Lemma~\ref{lem:oddeven}) uses a count of colors and the planarity of the board cut
open at the seam.

{\emergencystretch=5em
We begin with the colors. Color the cells as a checkerboard: the cell~$(r,c)$ has color~\mbox{$(r+c) \bmod 2$}. A domino
of $R$ covers one cell of each color. The seam domino~$s_b$ does not: both of its cells have color~$b \bmod 2$, since
$b + 2n \equiv \bar b + 1 \equiv b \pmod 2$. We use the coloring only through the
following count.\par}

\begin{lemma}\label{cl:par}
Every tiling of the M\"obius board uses as many seam dominoes~$s_b$ with $b$ even as
with $b$ odd. In particular, the number of seam dominoes it uses is even.
\end{lemma}

\begin{proof}
The seam dominoes~$s_b$ with $b$ even each cover two cells of color~$0$, those with
$b$ odd each cover two cells of color~$1$, and every other domino covers one cell of
each color. The board has equally many cells of each color, so a tiling uses equally
many seam dominoes of the two kinds.
\end{proof}

The count is similar to the one showing that a chessboard with two opposite corners
removed has no domino tiling. Here it reduces the question of types to a question about
the parities of the paired rows.

\begin{lemma}\label{lem:type}
If a seam cycle of $S$ pairs an even row with an odd row, then the two rows have the
same type.
\end{lemma}

\begin{proof}
Let $C$ be a seam cycle of $S$ that pairs rows~$b$ and $b'$, one even and one odd.
Suppose that one of the two rows has type~$h$ and the other has type~$s$. By
Proposition~\ref{thm:cars}, the only seam dominoes on $C$ are $s_b$ and $s_{b'}$. So
shifting along $C$ removes one of them, inserts the other, and changes no other seam
domino. Since one of $b$ and $b'$ is even and the other is odd, the shift replaces a
seam domino~$s_a$ with $a$ even by one with $a$ odd, or the reverse. By
Lemma~\ref{cl:par}, the tiling~$S$ uses as many seam dominoes of the first kind as of
the second. So $S \mathbin{\triangle} C$ uses two more of one kind than of the other,
and this contradicts Lemma~\ref{cl:par} for the tiling $S \mathbin{\triangle} C$.
\end{proof}

It remains to show that every seam cycle pairs an even row with an odd row. For this
we cut the board open at the seam. Every cell of the first column and every cell
of the last column is covered by exactly one of the seam dominoes $s_1, \dots, s_{2m}$.
We call these cells the \emph{border cells}, and we say that the two cells of $s_a$,
which are $(a,2n)$ and $(\bar a,1)$, have \emph{label}~$a$. The labels increase down
the last column and up the first column (Figure~\ref{fig:strands}).

Let $C$ be a seam cycle of $S$ that pairs rows~$b$ and $b'$. By
Proposition~\ref{thm:cars}, the only seam dominoes of $C$ are $s_b$ and $s_{b'}$.
Remove them from $C$. What is left is two disjoint paths, the \emph{strands} of $C$.
The dominoes of a strand are dominoes of $R$, so a strand is a \emph{path in} $R$: a sequence of cells of $R$ in which consecutive
cells share a side. Going round $C$, the two
seam dominoes separate the two strands. So each strand has one end at a cell of $s_b$
and the other end at a cell of $s_{b'}$. That is, each strand joins a border cell with
label~$b$ to a border cell with label~$b'$. In the running example, the seam cycle
pairs rows~$1$ and $2$. One of its strands joins $(1,4)$ to $(2,4)$, and the other
joins $(6,1)$ to $(5,1)$ (Figure~\ref{fig:running}, right).

We now show that the strands of $C$ separate the border cells whose labels lie strictly
between $b$ and $b'$ from the other border cells: a path in $R$ that contains no cell
of $C$ cannot join a cell of one kind to a cell of the other.

\begin{lemma}\label{lem:sep}
Let $C$ be a seam cycle of $S$ that pairs rows~$b$ and $b'$, where $b < b'$, and let
$Q$ be a path in $R$ that joins two border cells and contains no cell of $C$. If the
label of one end of $Q$ lies strictly between $b$ and $b'$, then so does the label of
the other end.
\end{lemma}

\begin{figrow}
\begin{minipage}[t]{\linewidth}\centering
\begin{tikzpicture}[x=8.2mm,y=8.2mm,line join=round,line cap=round]
\fill[crossodd] (6,-1.3) .. controls (4.3,-1.3) and (3,-1.7) .. (3.6,-2.3) .. controls (4.2,-2.9) and (4.4,-3.2) .. (6,-3.2) -- cycle;
\fill[crossodd] (0,-1.8) .. controls (1.5,-1.8) and (2.6,-2.3) .. (2,-2.8) .. controls (1.4,-3.3) and (1.2,-3.7) .. (0,-3.7) -- cycle;
\draw[border] (0,0) rectangle (6,-5);
\node[lab,anchor=east] at (-0.22,-0.3) {\scriptsize $2m$};
\node[lab,anchor=east] at (-0.22,-4.7) {\scriptsize $1$};
\node[lab,anchor=west] at (6.22,-0.3) {\scriptsize $1$};
\node[lab,anchor=west] at (6.22,-4.7) {\scriptsize $2m$};
\draw[pathA] (6,-1.3) .. controls (4.3,-1.3) and (3,-1.7) .. (3.6,-2.3) .. controls (4.2,-2.9) and (4.4,-3.2) .. (6,-3.2);
\draw[pathA] (0,-1.8) .. controls (1.5,-1.8) and (2.6,-2.3) .. (2,-2.8) .. controls (1.4,-3.3) and (1.2,-3.7) .. (0,-3.7);
\draw[pathB,densely dashed] (6,-1.95) .. controls (4.9,-1.95) and (4.9,-2.65) .. (6,-2.65);
\node[text=carlineB] at (4.75,-2.3) {\footnotesize $Q$};
\fill[pathAdot] (6,-1.3) circle (1.9pt);
\node[anchor=west,text=carlineA] at (6.22,-1.3) {\footnotesize $b$};
\fill[pathAdot] (6,-3.2) circle (1.9pt);
\node[anchor=west,text=carlineA] at (6.22,-3.2) {\footnotesize $b'$};
\fill[pathAdot] (0,-3.7) circle (1.9pt);
\node[anchor=east,text=carlineA] at (-0.22,-3.7) {\footnotesize $b$};
\fill[pathAdot] (0,-1.8) circle (1.9pt);
\node[anchor=east,text=carlineA] at (-0.22,-1.8) {\footnotesize $b'$};
\fill[pathBdot] (6,-1.95) circle (1.9pt);
\node[anchor=west,text=carlineB] at (6.22,-1.95) {\footnotesize $a$};
\fill[pathBdot] (6,-2.65) circle (1.9pt);
\node[anchor=west,text=carlineB] at (6.22,-2.65) {\footnotesize $a'$};
\fill[crossodd] (15.4,-1.3) .. controls (13.9,-1.3) and (13.6,-0.6) .. (12.6,-1) .. controls (11.6,-1.4) and (10.9,-1.8) .. (9.4,-1.8) -- (9.4,-3.7) .. controls (10.8,-3.7) and (11.4,-3.6) .. (12.4,-3.9) .. controls (13.4,-4.2) and (14,-3.2) .. (15.4,-3.2) -- cycle;
\draw[border] (9.4,0) rectangle (15.4,-5);
\node[lab,anchor=east] at (9.18,-0.3) {\scriptsize $2m$};
\node[lab,anchor=east] at (9.18,-4.7) {\scriptsize $1$};
\node[lab,anchor=west] at (15.62,-0.3) {\scriptsize $1$};
\node[lab,anchor=west] at (15.62,-4.7) {\scriptsize $2m$};
\draw[pathA] (15.4,-1.3) .. controls (13.9,-1.3) and (13.6,-0.6) .. (12.6,-1) .. controls (11.6,-1.4) and (10.9,-1.8) .. (9.4,-1.8);
\draw[pathA] (9.4,-3.7) .. controls (10.8,-3.7) and (11.4,-3.6) .. (12.4,-3.9) .. controls (13.4,-4.2) and (14,-3.2) .. (15.4,-3.2);
\draw[pathB,densely dashed] (15.4,-1.95) .. controls (14,-1.95) and (13,-3) .. (12,-2.5) .. controls (11,-2) and (10.6,-2.35) .. (9.4,-2.35);
\node[text=carlineB] at (12.7,-2.35) {\footnotesize $Q$};
\fill[pathAdot] (15.4,-1.3) circle (1.9pt);
\node[anchor=west,text=carlineA] at (15.62,-1.3) {\footnotesize $b$};
\fill[pathAdot] (15.4,-3.2) circle (1.9pt);
\node[anchor=west,text=carlineA] at (15.62,-3.2) {\footnotesize $b'$};
\fill[pathAdot] (9.4,-3.7) circle (1.9pt);
\node[anchor=east,text=carlineA] at (9.18,-3.7) {\footnotesize $b$};
\fill[pathAdot] (9.4,-1.8) circle (1.9pt);
\node[anchor=east,text=carlineA] at (9.18,-1.8) {\footnotesize $b'$};
\fill[pathBdot] (15.4,-1.95) circle (1.9pt);
\node[anchor=west,text=carlineB] at (15.62,-1.95) {\footnotesize $a$};
\fill[pathBdot] (9.4,-2.35) circle (1.9pt);
\node[anchor=east,text=carlineB] at (9.18,-2.35) {\footnotesize $a'$};
\end{tikzpicture}\\[6pt]
\begin{minipage}{.92\linewidth}\capt{\label{fig:strands}The two strands of a seam
cycle that pairs rows~$b$ and $b'$ are drawn in red on the rectangle~$R$, in the two
cases of the proof of Lemma~\ref{lem:sep}. The ends of the strands are marked with
their labels. The labels run from $1$ to $2m$ down the last column and up the first
column. The shaded region contains the border cells whose labels lie strictly between
$b$ and $b'$. The dashed blue curve is a path~$Q$ that joins two of these cells.}
\end{minipage}\end{minipage}
\end{figrow}

\begin{proof}
The rectangle~$R$ is a region of the plane. Draw each strand of $C$ as a curve through
the centers of its cells, continued at each end to the border of $R$, and draw the
path~$Q$ as a curve through the centers of its cells. Each curve stays inside the cells
of its path. Since $Q$ contains no cell of $C$, the curve of $Q$ does not meet
the curves of the strands. The curve of a strand does not cross itself, and it meets
the border of $R$ only at its two ends. So, by the Jordan curve theorem, it divides the
rest of $R$ into two parts, and $Q$ lies in one of them.

Each strand joins a border cell with label~$b$ to a border cell with label~$b'$. So
there are two cases (Figure~\ref{fig:strands}).

In the first case, one strand joins the two cells of the last column with labels~$b$
and $b'$, and the other strand joins the two cells of the first column with these
labels. In either column, the cells between the ends of the strand are the cells with
labels strictly between $b$ and $b'$. The strand divides those that are not on $C$
from every other border cell that is not on $C$.

In the second case, one strand joins the cell of the last column with label~$b$ to the
cell of the first column with label~$b'$, and the other strand joins the cell of the
last column with label~$b'$ to the cell of the first column with label~$b$. The labels
increase down the last column and up the first column. So in each of the two columns
the end of the first strand is above the end of the second, and the first strand is
the upper one. The border cells that have labels strictly between $b$ and $b'$ and are
not on $C$ are then the border cells that lie below the first strand and above the
second. Every other border
cell that is not on $C$ lies above the first strand or below the second.

In both cases, a path in $R$ that contains no cell of $C$ cannot join a border cell
whose label lies strictly between $b$ and $b'$ to a border cell whose label does not.
\end{proof}

We return to the tiling of Figure~\ref{fig:D}. Its red seam cycle pairs rows~$3$
and $10$. One of its strands joins $(3,20)$ to the cell~$(1,1)$, which has label~$10$,
and the other joins $(10,20)$ to the cell~$(8,1)$, which has label~$3$. So the red
cycle is in the second case of the proof. The shaded cells lie below its upper strand
and above its lower strand, and the shaded border cells are those with labels~$4$ to
$9$ that are not on the red cycle. One border cell with such a label is on the red
cycle: the cell~$(7,1)$, which has label~$4$. The blue cycle pairs rows~$6$ and $9$. Each of its strands is a path in $R$ that
contains no cell of the red cycle, and the blue cycle lies among the shaded cells, as
the lemma requires. The blue cycle itself is in the first case: each of its strands
returns to the column it left.

\begin{lemma}\label{lem:oddeven}
Every seam cycle pairs an even row with an odd row, and the two rows have the same type.
\end{lemma}

\begin{proof}
Let $C$ be a seam cycle of $S$ that pairs rows~$b$ and $b'$, where $b < b'$. We show
that the cells~$(a,2n)$ of the last column with $b < a < b'$ fall into pairs.

Let $(a,2n)$ be one of these cells. Suppose first that $a$ is not in $B(S)$. Then $S$
covers $(a,2n)$ by a vertical domino of the last column. The other cell of this domino
is in row~$a-1$ or row~$a+1$, and it is not in row~$b$ or row~$b'$, because these rows
are in $B(S)$. So its row is also strictly between $b$ and $b'$. We pair the two
cells of the domino.

Suppose instead that $a$ is in $B(S)$. Then a seam cycle~$C'$ of $S$ pairs row~$a$ with
another row~$a'$, and $C' \neq C$ because $a$ is neither $b$ nor $b'$. One strand of
$C'$ has an end at $(a,2n)$, which has label~$a$. This strand is a path in $R$ from
$(a,2n)$ to a border cell with label~$a'$, and it contains no cell of $C$, because
distinct seam cycles are disjoint. By Lemma~\ref{lem:sep}, the label~$a'$ lies strictly
between $b$ and $b'$. We pair $(a,2n)$ with $(a',2n)$.

There are $b' - b - 1$ cells~$(a,2n)$ with $b < a < b'$, and we have divided them into
pairs. So $b' - b - 1$ is even, and one of $b$ and $b'$ is even and the other is odd.
By Lemma~\ref{lem:type}, the two rows have the same type.
\end{proof}

We can now find the ordinary tilings.

\begin{proposition}[Shift classes]\label{thm:bij}
Every shift class contains exactly one ordinary tiling. If $T$ is an ordinary tiling,
then every marked domino of $T$ lies on a seam cycle of $T$, and every seam cycle of
$T$ contains exactly one marked domino.
\end{proposition}

\begin{proof}
Let $S$ be a M\"obius tiling. By Lemma~\ref{lem:oddeven}, each seam cycle of $S$ pairs
two seam rows or two rows of type~$h$. Shift along the seam cycles that pair seam
rows. These shifts reverse the type of every seam row and of no other row. So the
result~$T$ has no seam row, and $T$ is an ordinary tiling in the shift class of $S$.
Hence every shift class contains an ordinary tiling.

Suppose that a shift class contains ordinary tilings~$T$ and $T'$. The class is the
shift class of $T$, so $T' = T \mathbin{\triangle} \bigcup \mathcal{C}$ for a
set~$\mathcal{C}$ of seam cycles of $T$. If $\mathcal{C}$ contained a cycle, pairing
rows~$b$ and $b'$, then the shift along it would replace $h_b$ by $s_b$, and the shifts
along the other cycles of $\mathcal{C}$ would not change the type of row~$b$. So $T'$
would contain $s_b$. Hence $\mathcal{C}$ is empty, and $T' = T$.

Finally, let $T$ be an ordinary tiling. For each $b \in B(T)$, the tiling~$T$ contains
$h_b$ and its companion~$T^*$ contains $s_b$. These two dominoes share the
cell~$(b,2n)$, so they lie on the same component of $T \mathbin{\triangle} T^*$, which
is a seam cycle. In particular every marked domino of $T$ lies on a seam cycle. Now let
$C$ be a seam cycle of $T$. By Lemma~\ref{lem:oddeven} it pairs an even row~$b$ with an
odd row~$b'$. A marked domino covers a cell of the last column, and by
Lemma~\ref{lem:two} the dominoes of $T$ on $C$ that cover a cell of the last column are
$h_b$ and $h_{b'}$. Of these, only $h_b$ is marked.
\end{proof}

\pagebreak

We can now state the bijection.

\begin{theorem}[The bijection]\label{thm:phi}
For an ordinary tiling~$T$ and a marked domino~$d$ of $T$, let $C_d$ be the seam cycle
of $T$ that contains $d$. The map
\[
  \Phi(T,M) \;=\; T \mathbin{\triangle} \bigcup_{d \in M} C_d
\]
is a bijection from the set of pairs~$(T,M)$, where $T$ is an ordinary tiling and $M$
is a set of marked dominoes of $T$, to the set of M\"obius tilings.
\end{theorem}

\begin{proof}
Let $T$ be an ordinary tiling. By Proposition~\ref{thm:bij}, the assignment $d \mapsto C_d$
is a bijection from the marked dominoes of $T$ to the seam cycles of $T$. So as $M$
runs over the sets of marked dominoes of $T$, the tiling~$\Phi(T,M)$ runs over the
shift class of $T$, and it takes each tiling in the class once, because distinct sets
of seam cycles give distinct tilings. The shift classes partition the set of M\"obius
tilings, and by Proposition~\ref{thm:bij} each of them is the shift class of exactly one
ordinary tiling.
\end{proof}

\begin{corollary}\label{cor:main}
The number of tilings of the M\"obius board is $\sum_T 2^{k(T)}$, the sum running over
the tilings~$T$ of $R$.
\end{corollary}

\begin{proof}
An ordinary tiling~$T$ has $2^{k(T)}$ sets of marked dominoes.
\end{proof}

The tiling~$\Phi(T,M)$ is obtained from $T$ by pushing each domino of $M$ across the
seam and advancing the other dominoes on its seam cycle by one cell to make room, as
\cite{sandpiles} suggested. In the running example, $T$ has one marked domino, $h_2$,
and the shift class of $T$ consists of $T = \Phi(T,\emptyset)$ and
$T^* = \Phi(T,\{h_2\})$.

\section{Remarks}\label{sec:remarks}

The bijection evaluates neither side of \eqref{eq:main}. The right side is evaluated in
\cite{sandpiles}: for $n \ge 2$, Theorem~18 of \cite{sandpiles} shows that
$\sum_T 2^{k(T)}$ equals the product \eqref{eq:luwu}, by computing a determinant from
the theory of sandpiles by means of Chebyshev polynomials. That computation does not
use the Lu--Wu formula. In \cite{sandpiles} the Lu--Wu formula is needed only to
conclude that the product is the number of M\"obius tilings (Remark~20 of
\cite{sandpiles}). The bijection gives this conclusion directly. So the bijection and
the computation in \cite{sandpiles} together prove the Lu--Wu formula \eqref{eq:luwu}
for $n \ge 2$, as anticipated in \cite{sandpiles}. This proof uses no Pfaffians.

Code that enumerates the M\"obius tilings of a board, constructs companions and seam
cycles, and carries out the bijection is available
at \url{https://people.reed.edu/~davidp/moebius/}.
It confirms Propositions~\ref{thm:cars} and~\ref{thm:bij} and Theorem~\ref{thm:phi}
for every M\"obius tiling of the
boards it has been run on, the largest being the $8 \times 6$ board with $603{,}126$
tilings.

\medskip
\noindent\textbf{Declaration of generative AI use.}
The mathematics and the text of this article were worked out with Claude (Anthropic;
Claude Fable 5.2 and earlier models, through Claude Code) between August and November
2026.

The author
chose the problem and the kind of proof to look for, a bijection that moves dominoes as
\cite{sandpiles} suggested. Claude proposed the companion of a
tiling, which is the central construction of the article, and most of the arguments in
the proofs. Several times the
author judged a correct proof too complicated and asked for a simpler one. The
description of the companion as a reflection and the planar argument of
Section~\ref{sec:bij} came from these requests. Claude wrote the verification code and the code that generates the
figures, and it searched the literature for earlier work. The text was written
jointly, through ten complete drafts. Claude drafted, and the author read the drafts and directed
several hundred changes, on every scale from the organization of the article to the
wording of single sentences.

The author has checked every statement and proof, and takes full responsibility for the
content.

\bigskip
\noindent{\small\textsc{Department of Mathematics and Statistics, Reed College,
Portland, OR 97202, USA}\\
\textit{Email address:} \href{mailto:davidp@reed.edu}{\texttt{davidp@reed.edu}}}

\end{document}